\documentclass[11pt]{article}
\usepackage{latexsym}

\newtheorem{theorem}{Theorem}
\newtheorem{proposition}{Proposition}
\newtheorem{corollary}{Corollary}
\newtheorem{lemma}{Lemma}

\newtheorem{definition}{Definition}

\newtheorem{claim}{Claim}

\newtheorem{remark}{Remark}

\newcommand{\thmref}[1]{Theorem~\ref{thm:#1}} 
\newcommand{\lemref}[1]{Lemma~\ref{lem:#1}} 
\newcommand{\propref}[1]{Proposition~\ref{prop:#1}} 
\newcommand{\claimref}[1]{Claim~\ref{claim:#1}} 
\newcommand{\corref}[1]{Corollary~\ref{cor:#1}} 
\newcommand{\defref}[1]{Definition~\ref{def:#1}} 
\newcommand{\secref}[1]{Section~\ref{sec:#1}} 
\newcommand{\eqnref}[1]{(\ref{eq:#1})} 

\newcommand\ignore[1]{}

\newcommand\mypar[1]{\noindent\par{\sc #1}~}

\def\R{{\bf R}} 
\def\N{{\bf N}} 
\def\I{\sqrt{-1}}  

\def\periodeq{\mbox{.}}
\def\commaeq{\mbox{,\,}}

\renewcommand{\Pr}[1]{\mbox{\bf Pr}\left(#1\right)} 
\newcommand{\Ex}[1]{\mbox{\bf Ex}\left[#1\right]} 
\newcommand{\Prp}[2]{\mbox{\bf Pr}_{#1}\left(#2\right)} 
\def\eqdist{=^d}   

\newcommand{\bigoh}[1]{O\left(#1\right)}
\newcommand{\liloh}[1]{o\left(#1\right)}
\newcommand{\ohmega}[1]{\Omega\left(#1\right)}
\newcommand{\theita}[1]{\Theta\left(#1\right)}

\newcommand\QED{\ifhmode\allowbreak\else\nobreak\fi
\quad\nobreak$\Box$\medbreak}
\newcommand{\proofstart}{\par\noindent\sl Proof:\rm\enspace}
\newcommand{\proofend}{\QED\par}
\newenvironment{proof}{\proofstart}{\proofend}

\def\eps{\epsilon}

\def\mon{\mbox{\sf mon}}

\def\periodeq{\mbox{.}}
\def\commaeq{\mbox{,}}

\usepackage{latexsym}

\def\mon{\mbox{\sf Mon}}

\def\loserhasmorethan{\mbox{\sf LoserHasMoreThan}}
\def\hasmorethan{\mbox{\sf HasMoreThan}}

\begin{document}

\title{The Onset of Dominance in Balls-in-Bins Processes with Feedback}
\author{Roberto Imbuzeiro Oliveira\thanks{IBM T.J. Watson Research Center, Yorktown Heights, NY 10598. \texttt{riolivei@us.ibm.com}. Done while the author was a Ph.D. student at New York University under the supervision of Joel Spencer. Work funded a CNPq doctoral scholarship. Current address: IMPA, Rio de Janeiro Brazil.}} \maketitle
\begin{abstract}Consider a balls-in-bins process in which each new ball goes into a given bin with probability proportional to $f(n)$, where $n$ is the number of balls currently in the bin and $f$ is a fixed positive function. It is known that these so-called {\em balls-in-bins processes with feedback} have a monopolistic regime: if $f(x)=x^p$ for $p>1$, then there is a finite time after which one of the bins will receive all incoming balls.\\
Our goal in this paper is to quantify the onset of monopoly. We
show that the initial number of balls is large and bin $1$ starts
with a fraction $\alpha>1/2$ of the balls, then with very high
probability its share of the total number of balls never decreases
significantly below $\alpha$. Thus a bin that obtains more than
half of the balls at a ``large time" will most likely preserve its
position of leadership. However, the probability that the winning
bin has a non-negligible advantage after $n$ balls are in the
system is $\sim\mbox{const.}\times n^{1-p}$, and the number of
balls in the losing bin has a power-law tail. Similar results also
hold for more general functions $f$.\end{abstract}

\section{Introduction}\label{sec:intro}
Consider a discrete-time Markov process with $B$ bins, each one of
which containing $I_i(m)>0$ balls at time $m$ for each $m\in
\{0,1,2,\dots\}$ and $i\in\{1,\dots,B\}$. Its evolution is as
follows: at each time $m>0$, a ball is added to a bin $i_m$, so
that $I_{i_m}(m)=I_{i_m}(m-1)+1$ and $I_j(m)=I_j(m-1)$ for all
$i\in\{1,\dots,B\}\backslash\{i_m\}$, and the random choice of bin
$i_m$ has distribution
\begin{equation}\label{eq:recipe}\Pr{i_{m}=i\mid \{I_j(m-1)\,:\,{1\leq j\leq B}\}}=
\frac{f(I_i(m-1))}{\sum_{j=1}^B f(I_j(m-1))}\;\; (1\leq i\leq
B),\end{equation} for a fixed positive function $f:\N\to
(0,+\infty)$. This recipe specifies what we call a {\em
balls-in-bins process with feedback function $f$ and $B$ bins},
and one should notice that when $f$ is an increasing function
--- the case that has been mostly considered in the literature ---
there is a tendency that {\em the rich get richer}: the more balls
a bin has, the more likely it is to receive the next ball.

This class of processes\footnote{A longer background discussion is
available from the first author's PhD thesis \cite{tese}.} was
proposed by Drinea, Frieze and Mitzenmacher \cite{Drinea02} as a
model for competing products in an economy, as well as a simpler
variant of so-called preferential-attachment models for large
networks (see \cite{AlbertSurvey} for a survey of the latter). That
paper concentrates on the special case where $f(x)=x^p$ for some
parameter $p>0$. The authors prove that when $p>1$, there almost
surely exists one bin that gets all but a negligible fraction of the
balls in the large-time limit; whereas for $p<1$, the asymptotic
fractions of balls in each bin are all the same. The $p=1$ case is
the classic P\'olya Urn model, for which it has been long known that
the number of balls in each bin converges almost surely to a
non-degenerate random variable, and thus the process has different
{\em regimes} depending on whether $p<1$, $p=1$ or $p>1$.

However, stronger results on the $p>1$ case are available. A paper
by Khanin and Khanin \cite{Khanin01} introduced what amounts to
the same process as a model for neuron growth, and proved the
following stronger result: if $p>1$, there almost surely is some
bin that gets {\em all but finitely many balls}. That is to say,
consider the following event, in which bin $i$ is the only one to
receive balls after some finite time $M$ (we call this {\em
monopoly by bin $i$}).
\begin{eqnarray}\mon_i &\equiv& \{\exists M\in \N\;\forall m\geq M\;\;\forall j\in [B]\;\;j\neq i\Rightarrow I_j(m)=I_j(M)\}\\ &=&\{\exists M\in \N\;\forall m\geq M\;\;i_{m}=i\}.\end{eqnarray}
The result of \cite{Khanin01} -- or rather, a straightforward
extension of it proven in \cite{Spencer??,tese} -- says that
\begin{theorem}[From \cite{Khanin01,Spencer??,tese}]\label{thm:Khanin}If
$\{I_m\}_{m=0}^{+\infty}$ is a balls-in-bins process with $B$ bins
and feedback function $f=f(x)>0$ satisfying
\begin{equation}\label{eq:monopoly}\sum_{j=1}^{+\infty}
\frac{1}{f(j)} <+\infty.\end{equation}satisfies $$\Pr{\exists
i\in[B]\,:\,\mon_i}=1.$$ In particular, this holds for $f(x)=x^p$,
$p>1$.\end{theorem} This is a much stronger statement than the one
that Drinea et al. proved in \cite{Drinea02}, and we take it as
our starting point in the present paper. We sketch a proof of
\thmref{Khanin} in \secref{exp_monopoly} below, since it helps to
build some intuition for our own results.

Our main goal is to quantify the monopoly phenomenon in the case
where there are $B=2$ bins. As powerful as \thmref{Khanin} is, it
tells us nothing about {\em how fast} the system approaches this
asymptotic regime. Informally, we will be interested in questions
like: suppose we throw in a million balls into $B=2$ bins. Should
we expect to find a bin with $10 \%$ more balls than the other
bin? And in case that does happen, is the leading bin likely to
lose its lead as we add more and more balls into the system? This
paper treats rigorous forms of those problems, for a broad (but
not entirely general) class of feedback functions. A summary of
our results is given below.

\begin{enumerate}
\item{\em Imbalanced start.} Start the process with a total of
$t\gg 1$ balls in the two bins and at least $\alpha=52\%$ of the
balls in bin $2$. We show in \secref{tails_imbalanced} that with
very high probability, there is no future time at which bin $2$
will have $\beta=51\%$ of the total number of balls in this
balls-in-bins process. A similar result holds for any other
$\alpha>\beta>1/2$. \ignore{More generally, if the initial number
of balls is large, and one of the bins has a fraction $1/2+\eps$
of all of the balls at time $0$, not only will that bin most
likely achieve monopoly, but also it will retain its advantage
over the other bin at all times, with very high probability.}
\item{\em Balls in the losing bin.} Let $L$ be the number of balls
that go into the {\em losing} bin, i.e. the bin that does not
achieve monopoly, when the initial number of balls is fixed. We
prove in \secref{tails_losing} that the distribution of $L$ has a
{\em heavy tail}. More concretely, if $f(x)=x^p$ for some $p>1$,
$\Pr{L>n}\sim c_p\times n^{1-p}$ for large $n$. \item{\em The time
until imbalance.} We show in \secref{tails_fraction} that for
fixed initial conditions, the probability that the losing bin has
at least $\alpha n$ balls at time $n$ (for $\alpha<1/2$ fixed)
also decays slowly in $n$. In particular, in the case $f(x)=x^p$,
this probability is $\sim c'_p\times n^{1-p}$ for $n$ large; that
is, we have a power law with the same exponent as in $2.$. An
extension of this result is presented in \secref{last}. \ignore{We
also show that the probability that the losing bin has $\geq
n/2-\lambda\sqrt{n}$ balls for a constant $\lambda$ decays like
$n^{1/2-p}$.}
\end{enumerate}

The picture that emerges from these results is that {\em it takes
a long time before a clear leader emerges, but once it does, it is
likely to stick}. One indication of the heavy-tail part of our
results was in the numerical simulations of \cite{Drinea02}, which
indicated that the a clear leader of the process took a long time
to emerge. Reference \cite{Drinea02} also showed that once the
leader emerges, it achieves dominance very quickly, but our
results regrading this (cf. $1.$) are stronger. \ignore{We also
note that the ``square-root separation" alluded to in item $3.$ is
not at all artificial. It is shown in \cite{MitzenmacherOS04} that
if each bin starts out with $n/2\pm \lambda\sqrt{n}$ balls, $n\gg
1$, there is a non-trivial probability of monopoly going either
way.}

Our main technical tool has been employed in
\cite{Khanin01,Spencer??} and other works, and seems to have
originated in Davis' work on reinforced random walks
\cite{Davis90}. We shall {\em embed} the discrete-time process we
are interested in into a continuous-time process built from
exponentially distributed random variables. The most salient
feature of this so-called {\em exponential embedding} is that
arrival times at different bins are independent and have an
explicit distribution. This greatly simplifies calculations and
permits the use of Chernoff-like bounds that we develop below.

The remainder of the paper is organized as follows. We discuss
preliminary material in \secref{prelim}. \secref{exp_embed}
rigorously introduces the exponential embedding process and
discusses its key properties (including \thmref{Khanin}). In
\secref{technical} we detail the assumptions we make on our
feedback functions $f$, while also deriving some consequences of
those assumptions. The next three sections correspond to items
$1.$---$3.$ . \secref{last} discusses possible extensions to our
results and some related work.

\mypar{Acknowledgements.} I thank Michael Mitzenmacher and Eleni
Drinea at Harvard University for early discussions regarding this
work. Thanks also go to my former Ph.D. advisor Joel Spencer, who
introduced me to this topic and made numerous valuable suggestions
to the work presented here.

\section{Preliminaries}\label{sec:prelim}

\par {\em Set notation.} Throughout the paper,
$\N=\{1,2,3,\dots\}$ is the set of positive integers, $\R^+ =
[0,+\infty)$ is the set of non-negative reals, and for any
$k\in\N$ $[k]=\{1,\dots,k\}$.\\

\par {\em Asymptotics.} We will use the standard $\bigoh{\cdot}/\liloh{\cdot}/\ohmega{\cdot}/\ll/\sim$ asymptotic notation. Let $f,g$ be functions of a real parameter $t$ and $t_0$ be a limit point of the domain of the two functions. We will say that $f(t)=\bigoh{g(t)}$ (or equivalently $g(t)=\ohmega{f(t)}$) as $t\to t_0$ when $\limsup_{t\to t_0} |f(t)/g(t)|<+\infty$. We will also say that $f(t)=\liloh{g(t)}$ (or $f(t)\ll g(t)$) as $t\to t_0$ when $\limsup_{t\to t_0} |f(t)/g(t)|=0$. Finally, we will say that $f(t)\sim g(t)$ as $t\to t_0$ when $\lim_{t\to t_0}f(t)/g(t)=1$.\\

\par {\em Balls-in-bins.} Formally, a balls-in-bins process with feedback function $f:\N\to (0+\infty)$ and $B\in\N$ bins is a discrete-time Markov chain $\{(I_1(m),\dots,I_B(m))\}_{m=0}^{+\infty}$ with state space $\N^B$ and transition probabilities given in the Introduction (see \eqnref{recipe}). We will usually refer to the index $i_m\in[B]$ as {\em the bin that receives a ball at time $m$}.

If $B=2$, $n\in\N$ and $0\leq \alpha \leq 1$, we sometimes denote by
$[n,\alpha]$ the state $(\lceil \alpha n\rceil,n-\lceil \alpha
n\rceil)\in\N^2$, i.e. there are $n$ bins in the system and an
$\alpha$-fraction of them is in bin $1$ (with rounding). This
alternative notation will be used whenever convenient.

For any $B$, if $E$ is an event of the process and $u\in\N^B$,
$\Prp{u}{E}$ is the probability of $E$ when the initial conditions
are set to $u$. The same notation will be used for the exponential
embedding defined in \secref{exp_embed}.\\

\par {\em Exponential random variables.} $X\eqdist \exp(\lambda)$
means that $X$ is a random variable with exponential distribution
with rate $\lambda>0$, meaning that $X\geq 0$ and
$$\Pr{X>t} = e^{-\lambda t}\;\; (t\geq 0).$$
The shorthand $\exp(\lambda)$ will also denote a generic random
variable with that distribution. Some elementary but extremely
useful properties of those random variables include
\begin{enumerate}\item {\em Lack of memory.} Let $X\eqdist \exp(\lambda)$
and $Z\geq 0$ be independent from $X$. The distribution of $X-Z$
conditioned on $X>Z$ is still equal to $\exp(\lambda)$. \item {\em
Minimum property.} Let $\{X_i\eqdist \exp(\lambda_i)\}_{i=1}^{m}$
be independent. Then $$X_{\min}\equiv \min_{1\leq i\leq
m}X_i\eqdist \exp(\lambda_1+\lambda_2+\dots \lambda_m)$$ and for
all $1\leq i\leq m$
\begin{equation}\Pr{X_i = X_{\min}} = \frac{\lambda_i}{\lambda_1+\lambda_2+\dots
\lambda_m}\end{equation}  \item {\em Multiplication property.} If
$X\eqdist\exp(\lambda)$ and $\eta>0$ is a fixed number, $\eta X
\eqdist\exp(\lambda/\eta)$. \item {\em Moments and transforms.} If
$X\eqdist\exp(\lambda)$, $r\in\N$ and $t\in\R$,
\begin{eqnarray}
\Ex{X^r} & = & \frac{r!}{\lambda^r}\commaeq \\
\Ex{e^{\I t X}} & = & \frac{1}{1 - \frac{\I t}{\lambda}}\\
\Ex{e^{t X}} & = & \left\{\begin{array}{ll}\frac{1}{1 - \frac{t}{\lambda}} & (t<\lambda)\\
                    +\infty  & (t\geq
                    \lambda)\end{array}\right.\end{eqnarray}

\end{enumerate}

\section{The exponential embedding}\label{sec:exp_embed}

\subsection{Definition and key properties}

Let $f:\N\to(0,+\infty)$ be a function, $B\in\N$ and
$(a_1,\dots,a_B)\in\N^B$. We define below a continuous-time
process with state space $(\N\cup\{+\infty\})^B$ and initial state
$(a_1,\dots,a_B)$ as follows. Consider a set $\{X(i,j)\,:\,i\in
[B],\, j\in \N\}$ of independent random variables, with
$X(i,j)\eqdist \exp(f(j))$ for all $(i,j)\in [B]\times \N$, and
define

\begin{equation}\label{eq:exp_contproc}N_i(t)\equiv \sup\left\{n\in
\N\,:\,\sum_{j=a_i}^{n-1} X(i,j) \leq t\right\}\;\;\;(i\in
[B],t\in\R^+ = [0,+\infty))\commaeq\end{equation} where by
definition $\sum_{j=i}^{k}(\dots)=0$ if $i>k$. Thus $N_i(0)=a_i$
for each $i\in [B]$, and one could well have $N_i(T)=+\infty$ for
some finite time $T$ (indeed, that {\em will} happen for our cases
of interest); but in any case, the above defines a continuous-time
stochastic process, and in fact the $\{N_i(\cdot)\}_{i=1}^B$
processes are independent. Each one of this processes is said to
correspond to {\em bin} $i$, and each one of the times
$$X(i,a_i),X(i,a_{i})+X(i,a_{i}+1),X(i,a_{i})+X(i,a_{i}+1)+X(i,a_{i}+2),\dots$$
is said to be an {\em arrival time at bin $i$}. As in the
balls-in-bins process, we imagine that each arrival correspond to
a ball being placed in bin $i$.

In fact, we {\em claim} that this process is related as follows to
the balls-in-bins process with feedback function $f$, $B$ bins and
initial conditions $(a_1,\dots,a_B)$.

\begin{theorem}[Proven in \cite{Davis90,Khanin01,Spencer??,tese,Oliveira07Brown}]\label{thm:exp_embed}Let the $\{N_i(\cdot)\}_{i\in[B]}$ process be defined as above. One can order the arrival times of the $B$ bins in increasing order (up to their first accumulation point, if they do accumulate) so that $T_1<T_2<\dots$ is the resulting sequence. The distribution of
$$\{I_m = (N_1(T_m),N_2(T_m),\dots,N_B(T_m))\}_{m\in\N}$$
is the same as that of a balls-in-bins process with feedback
function $f$ and initial conditions
$(a_1,a_2,\dots,a_B)$.\end{theorem}

One can prove this result\footnote{The exact attribution of this
result is somewhat confusing. Ref. \cite{Khanin01} cites the work
of Davis \cite{Davis90} on reinforced random walks, where it is in
turn attributed to Rubin.} as follows. First, notice that the {\em
first arrival time $T_1$} is the minimum of $X(j,a_j)$, ($1\leq
j\leq B$). By the minimum property presented above, the
probability that bin $i$ is the one at which the arrival happens
is like the first arrival probability in the corresponding
balls-in-bins process with feedback:
\begin{equation}\label{eq:min}\Pr{X(i,a_i) = \min_{1\leq j\leq
B}X(j,a_j)} = \frac{f(a_i)}{\sum_{j=1}^B f(a_j)}.\end{equation}
More generally, let $t\in\R^+$ and condition on $(N_i(t))_{i=1}^B
= (b_i)_{i=1}^B \in\N^B$, with $b_i\geq a_i$ for each $i$ (in
which case the process has not blown up). This amounts to
conditioning on
$$\forall i\in[B]\;\; \sum_{j=a_i}^{b_i-1}X(i,b_i)\leq t < \sum_{j=a_i}^{b_i}X(i,b_i).$$
From the lack of memory property of exponentials, one can deduce
that the first arrival after time $t$ at a given bin $i$ will happen
at a $\exp(f(b_j))$-distributed time, independently for different
bins. This takes us back to the situation of \eqnref{min}, only with
$b_i$ replacing $a_i$, and we can similarly deduce that bin $i$ gets
the next ball with the desired probability,
$$\frac{f(b_i)}{\sum_{j=1}^Bf(b_j)}.$$

\subsection{The occurrence of monopoly}\label{sec:exp_monopoly}

The exponential embedding yields a ``Book proof" of
\thmref{Khanin}. Without loss of generality, we can assume that
$B=2$ (the general case follows from comparing pairs of bins). The
notation we use comes from the previous section, and we also
employ the version of the balls-in-bins process given by
$I_m=(N_i(T_m))_{i=1,2}$ (cf. \thmref{exp_embed}).

Under the condition that
\begin{equation}\label{eq:monopoly_cond}\sum_{j=1}^{+\infty}
\frac{1}{f(j)} <+\infty,\end{equation} one has
\begin{equation} F_i\equiv \sum_{j=a_i}^{+\infty}X(i,j)<+\infty\mbox{ almost surely  } (i=1,2).\end{equation}
Indeed, the terms in $F_i$ are all positive and since
$X(i,j)\eqdist \exp(f(j))$ for all $i,j$,
$$\Ex{F_i} = \sum_{j=a_i}^{+\infty} \frac{1}{f(j)}<+\infty\mbox{  by \eqnref{monopoly_cond}.}$$
It is also easy to see that the $F_i$'s  are independent and have
no point masses in their distributions. Thus with probability $1$,
either $F_1<F_2$ or $F_2<F_1$.

Suppose that the first alternative holds. Since
$\sum_{j=a_i}^{N-1} X(i,j)\nearrow F_i$ as $N\to +\infty$, we can
deduce that

\begin{equation}\exists n_2\in\N \;\;\;\forall n\geq a_1 \;\; A_{1,n} \equiv\sum_{j=a_1}^{n-1}X(1,j)<F_1<A_{2,n_2}\equiv \sum_{\ell=a_2}^{n_2-1}X(2,j).\end{equation}

The sequence $\{A_{1,n}\}_{n\in\N}$ is composed of arrival times;
that is to say, it is a subsequence of $\{T_m\}_{m\in\N}$.
Moreover, that sequence converges to $F_1$. It follows that the
first accumulation point of the sequence $\{T_m\}_{m\in\N}$ is at
most $F_1$, and that $T_m\leq F_1$ for all $m$. But since
$A_{2,n_2}>F_1$, this implies that for all times $m$ of the
discrete-time process,
$$I_m(2)=N_2(T_m)<n_2\periodeq$$ Thus if $F_1<F_2$, there exists a finite
$n_2\in\N$ such that $I_m(2)<n_2$ for all $m$, i.e. bin $2$ never
has more than $n_2$ balls. This implies that bin $1$ {\em must}
achieve monopoly whenever $F_1<F_2$.

If $F_1>F_2$, the same reasoning shows that bin $2$ achieves
monopoly. As pointed out above, with probability $1$ either
$F_1<F_2$ or $F_2>F_1$; thus the proof is finished.

\begin{remark}It is not hard to show that $F_i=+\infty$ almost surely if
$\sum_j f(j)^{-1}=+\infty$, and in that case monopoly has
probability $0$. The interested reader can see
\cite{Khanin01,tese,Oliveira07Brown} for details.\end{remark}

\begin{remark}\label{rem:fictitious}Assume that bin $1$ achieves monopoly, as in the proof above. In this case, all arrivals of the {\em continuous-time} process at bin $2$ after time $F_1$ do not actually happen in the embedded {\em discrete-time} process $\{I_m=(I_m(1),I_m(2))\}$. We call these ``ghost events" a {\em fictitious continuation} of our process. This very useful device is akin to the continuation of a Galton-Watson process beyond its extinction time (see e.g. \cite{AlonSpencer_Method}) and is equally useful in calculations and proofs.\end{remark}

\section{Assumptions on feedback functions and a large-deviations bound}\label{sec:technical}

The purpose of this rather technical section is two-fold. First,
we spell out the technical assumptions on the feedback function
$f$ that we need in our proofs. Nothing seems to actually {\em
require} these assumptions, but they facilitate certain estimates
that we employ in the proofs.

The second purpose is to present a large-deviations bound on
random sums such as $\sum_{j=m}^{+\infty}X(i,j)$. One can check
that the {\em variance} of such a sum decreases as $m\to +\infty$;
we show below that these sums are in fact close to their means
with all-but-exponential probability.

Some readers might wish to skip the proofs in this section on a
first reading.

\subsection{Valid feedback functions}

The feedback functions we allow in our results satisfy the
following definition.

\begin{definition}\label{def:est_smooth} An increasing function $f:\N\to (0,+\infty)$ with $f(1)=1$\footnote{The requirement that $f(1)=1$ is just a normalization condition,
as it does not change the process.} is said to be a {\em valid
feedback function} if it can be extended to a $C^{1}$ function
$g:[1,+\infty)\to(0,+\infty)$ with the following property: if $(\ln
g(\cdot))'$ is the right-derivative of $\ln g$, and $h(x)\equiv
x(\ln g(x))'$ (for $x\in\R^+\cup\{0\}$),
    \begin{enumerate}
        \item $\liminf_{x\to +\infty}h(x)>1$;
        \item $\lim_{x\to +\infty}x^{-{1/4}}h(x) = 0$;
        \item there exist $C>0$ and $x_0\in\R^+$ such that for all $\eps\in(0,1)$ and all $x\geq x_0$
\begin{equation}\label{eq:est_hslow}\sup_{x\leq t\leq
x^{1+\eps}}\left|\frac{h(t)}{h(x)}-1\right|\leq
C\,\eps\periodeq\end{equation}
    \end{enumerate}
With slight abuse of notation, we will always assume that $f$ is
defined over $[1,+\infty)$ and is $C^1$. We will also call $h$ the
{\em characteristic exponent} of $f$.\end{definition} Functions with
exponential growth (such as $f(x)=2^x$) or with oscillations fail to
satisfy
\defref{est_smooth}. On the other
hand, requiring that $f$ be increasing seems natural, and the
assumption still leaves us with plenty of interesting examples of
feedback functions. For instance, any of the functions defined below
\begin{eqnarray*}f(x) &=&x^p \;(\mbox{ for some fixed }p>1)\commaeq\\
f(x) & =& x^{p\ln^\alpha x} \;(\mbox{ for some fixed }p>1,\alpha>0)\commaeq\\
f(x) &=& x^{p}\ln(x+e-1)\;(\mbox{ for some fixed
}p>1)\periodeq\end{eqnarray*} is valid. The ``canonical case'' where
$f(x)=x^p$ ($x\geq 1$) explains the terminology for the
characteristic exponent: in that case, $h(x)\equiv p$ for all $x>1$.
We also note that whenever $f$ is a valid feedback function, the
monopoly condition is satisfied.

\begin{proposition}\label{prop:est_mon}If $f$ is a valid feedback function,
$\sum_{j\in\N}f(j)^{-1}<+\infty$.\end{proposition}
\begin{proof}The condition $\liminf_{x\to +\infty}h(x)>1$ implies that there exists a $n\in\N$ such that $h(x)\geq
c>1$ for all $x\geq n$. This implies that $f(j)=\ohmega{j^c}$ as
$j\to +\infty$, which is enough for the convergence of
$\sum_{j\in\N}f(j)^{-1}$.\end{proof}

\subsection{Consequences of the definition}

Let us now define the quantity
\begin{equation}S_r(n,m) \equiv \sum_{j=n}^{m-1}\frac{1}{f(n)^r}\;\;\;(r\in\R^+,\,n\in\N, m\in\N\cup\{+\infty\})\end{equation}
for some $f:\N\to (0,+\infty)$, and also let $S_r(n)\equiv
S_r(n,+\infty)$ (which might diverge for some $r$). If $f(x)=x^p$
and $r\geq 1$, a simple calculation shows that for $n\gg 1$
$$S_r(n)\sim\int_{n}^{+\infty} \frac{dx}{f(x)^r} = \frac{n^{1-rp}-m^{1-rp}}{(rp-1)}.$$
The main content of the following lemma is that a similar result
holds for any valid $f$, if $p$ is replaced by the characteristic
exponent $h$.
\begin{lemma}\label{lem:est_moment}Assume that $f$ is a valid feedback function with characteristic exponent $h$.
Define the possibly divergent integrals:
$$M_r(n) = \int_n^{+\infty} \frac{dx}{f(x)^r}\;\;\;\;\;(r\in\R^+,n\in\N)\periodeq$$
Then for all $r\geq 1$ both $S_r(n)$ and $M_r(n)$ converge and
moreover, as $n\to +\infty$
$$S_r(n) \sim M_r(n) \sim \frac{n}{(rh(n)-1)f(n)^r}\periodeq$$
and for all fixed $r\geq 1,\rho>1$ there exists $a<1$ such that
$$M_r(\rho x)\leq aM_r(x)\mbox{ for all large $x$}\periodeq$$
Thus for any $r\geq 1$ and $\rho>1$:
$$S_r(n,\lceil\rho n\rceil)=\ohmega{\frac{n}{(rh(n)-1)f(n)}}.$$
\end{lemma}

Before we present the proof of \lemref{est_moment}, we state two
other lemmas. They follow directly from the assumptions on $h(x)$
and we omit their proofs.

\begin{lemma}\label{lem:est_bigger}$S_1(n)\gg e^{-n^{1/4}}$ as $n\to +\infty$.\end{lemma}
\begin{lemma}\label{lem:est_bounded}For any bounded function $w:\N\to\R$, $f(n + w(n))\sim f(n)$ and $h(n+w(n))\sim h(n)$ as $n\to +\infty$.\end{lemma}

\begin{proof}[of \lemref{est_moment}] Essentially the same proof appears in the last result in \cite{MitzenmacherOS04}, but we reproduce the argument here for convenience. Under our assumptions, $F(\cdot)\equiv f(\cdot)^r$ is a valid feedback function with characteristic exponent $rh(\cdot)$. Clearly, the lemma holds for $f$ iff it holds for $F$ with $r=1$. It follows that it suffices to prove the lemma in the case $r=1$, which is what we do below.

We can assume that $n$ is large enough, so that $\inf_{x\geq
n}h(x)=c>1$, and moreover
\begin{equation}\label{eq:est_hslow2}\forall x\geq n\commaeq \forall \eps>0 \;\; \sup_{x\leq t\leq x^{1+\eps}}\left|\frac{h(t)}{h(x)}-1\right|\leq C\eps,\end{equation}
which follows from the assumptions on $h$ in \defref{est_smooth}.
We start by noting that for all $y>x\geq n$
\begin{equation}\label{eq:est_fracgrowth}\frac{f(y)}{f(x)}\geq \left(\frac{y}{x}\right)^{c}\end{equation}
To see this, it suffices to notice that
$$\ln\frac{f(y)}{f(x)} = \int_{x}^y (\ln f(u))' \, du = \int_{x}^y \frac{h(u)}{u} \, du\geq (\inf_{x\leq u\leq y}h(u))\,\ln\frac{y}{x} \geq c\ln\frac{y}{x}\periodeq$$
Inequality \eqnref{est_fracgrowth} implies that $f(x)\gg x$ as
$x\to +\infty$. This justifies the following integration by parts
procedure.
\begin{eqnarray}\label{eq:est_firstI}M_1(n) &=& \int_{n}^{+\infty}\frac{dx}{\exp(\ln f(x))} \\
 & = & \left.\frac{x}{f(x)}\right|_{x=n}^{x\to +\infty} +  \int_{n}^{+\infty}\frac{x (\ln f(x))'\, dx}{\exp(\ln f(x))} \\ \label{eq:est_lastI}
& = & - \frac{n}{f(n)} +
\int_{n}^{+\infty}\frac{h(x)\,dx}{f(x)}\end{eqnarray} where in the
last line we plugged in the definition of $h$. We now {\em claim}
that
\begin{equation}\label{eq:est_claimhout}\int_{n}^{+\infty}\frac{h(x)\,dx}{f(x)} \sim  h(n)\,\int_{n}^{+\infty}\frac{dx}{f(x)} = h(n)\,M_1(n) \mbox{ as }n\to +\infty\periodeq\end{equation}
We will prove \eqnref{est_claimhout} below, but first we show how
it implies the lemma. Employing \eqnref{est_claimhout} with
equations \eqnref{est_firstI} to \eqnref{est_lastI} shows that
\begin{equation}\label{eq:est_T1OK}M_1(n) \sim \frac{n}{(h(n)-1)f(n)}\mbox{ as }n\to +\infty\commaeq\end{equation}
since $h(n)\geq c>1$ for $n$ large (as discussed above). By the
smoothness assumption, $h(n)/n\to 0$ as $n\to +\infty$, and
therefore
\begin{equation}\label{eq:est_T1bigger}\frac{n}{(h(n)-1)f(n)}\gg \frac{1}{f(n)}\mbox{ as }n\to +\infty\periodeq\end{equation}
Now notice that, since $f$ is increasing,
$$- \frac{1}{f(n)}\leq M_1(n) - S_1(n) \leq 0\commaeq$$ hence $|M_1(n)-S_1(n)|\leq
f(n)^{-1}$ and by equations \eqnref{est_T1OK} and
\eqnref{est_T1bigger},
$$S_1(n)\sim M_1(n) \sim \frac{n}{(h(n)-1)f(n)}\mbox{ as }n\to +\infty\commaeq$$
as desired. Moreover, if $\rho>1$ is fixed, equations
\eqnref{est_hslow2} and \eqnref{est_fracgrowth} imply that, for
large $n$,
$$h(\lceil\rho n\rceil)\sim h(n)\mbox{ and }f(\lceil\rho n\rceil)\geq \rho^c f(n)\mbox{ with }c \mbox{ as above}\commaeq$$
hence
$$\frac{\rho n}{(h(\lceil\rho n\rceil)-1)f(\lceil\rho n\rceil)} \leq (1+\liloh{1})\rho^{1-c}\frac{n}{(h(n)-1)f(n)}\commaeq$$
from which the desired statement about $M_r(\lceil\rho n\rceil)$
follows.

We now prove \eqnref{est_claimhout}. Choosing an arbitrary (but
fixed) $\eps>0$, we first show that
\begin{equation}\label{eq:est_smalltail}\mbox{as }n\to +\infty\commaeq \int_{n^{1+\eps}}^{+\infty}\frac{h(x)\,dx}{f(x)}\ll \int_{n}^{+\infty}\frac{h(x)\,dx}{f(x)}\periodeq\end{equation}
Indeed,
\begin{eqnarray*}\int_{n^{1+\eps}}^{+\infty}\frac{h(x)\,dx}{f(x)} &=& (1+\eps)\int_{n}^{+\infty}\frac{h(u^{1+\eps})\,u^\eps\,du}{f(u^{1+\eps})}\\
&\leq &(1+\eps)(1+C\eps)\int_{n}^{+\infty}\frac{h(u)\,u^\eps\,du}{f(u^{1+\eps})}\\
&\leq &(1+\eps)(1+C\eps)\int_{n}^{+\infty}\frac{h(u)\,u^\eps\,du}{f(u)u^{c\eps}}\\
& \leq & (1+\eps)(1+C\eps)n^{(1-c)\eps}\,\int_{n}^{+\infty}\frac{h(u)\,du}{f(u)}\\
&\ll &
\int_{n}^{+\infty}\frac{h(u)\,du}{f(u)}\commaeq\end{eqnarray*} where
the first line is a change of variables, the second line employs
\eqnref{est_hslow2}, the third line uses \eqnref{est_fracgrowth}
applied to $x=u$ and $y=u^{1+\eps}$, and the remaining lines follow
from $c>1$. The sequence of equations proves \eqnref{est_smalltail},
which implies in particular that, for $n$ large,
$$\int_{n}^{+\infty}\frac{h(x)\,dx}{f(x)}\sim \int_{n}^{n^{1+\eps}}\frac{h(x)\,dx}{f(x)} $$
But another use of \eqnref{est_hslow2} implies that $$(1 -
C\eps)h(n)\leq
\frac{\int_{n}^{n^{1+\eps}}\frac{h(x)\,dx}{f(x)}}{\int_{n}^{n^{1+\eps}}\frac{dx}{f(x)}}
\leq (1 + C\eps)h(n)\commaeq$$ and, similarly to
\eqnref{est_smalltail}, one can show that $$\mbox{as }n\to
+\infty\commaeq \int_{n^{1+\eps}}^{+\infty}\frac{dx}{f(x)}\ll
\int_{n}^{+\infty}\frac{dx}{f(x)}\periodeq$$ Thus we conclude that,
as $n\to +\infty$, $$1 - C\eps - \liloh{1}\leq
\frac{\int_{n}^{+\infty}\frac{h(x)\,dx}{f(x)}}{h(n)\int_{n}^{+\infty}\frac{dx}{f(x)}}\leq
1 + C\eps + \liloh{1}\periodeq$$ Since $\eps>0$ is arbitrary,
\eqnref{est_claimhout} follows, and the proof is
finished.\end{proof}

\subsection{A large-deviations estimate}

If bin $1$ starts with $a_1$ balls, the time until it has $b_1>a_1$
balls in the continuous-time process is $$\sum_{j=a_1}^{b_1-1}
X(1,j)$$ and the time until bin $1$ acquires infinitely many balls
is $$\nonumber\sum_{j=a_1}^{+\infty} X(1,j).$$ The latter is a sum
of independent $\exp(f(j))$ random variables, and it converges
whenever $\sum_{j} f(j)^{-1}<+\infty$ (cf. \secref{exp_monopoly}).
We will show that the sum concentrates very strongly around its
mean.

\begin{lemma}\label{lem:exp_lgdev}Let $f$ be a valid feedback function, so that the monopoly condition $\sum_jf(j)^{-1}<+\infty$ holds (cf. \propref{est_mon}). Assume that $\{V_j\eqdist\exp(f(j))\}_{j\in\N}$ be
a sequence of independent random variables, and define (for
$n\in\N$)
$$A_{n} \equiv \sum_{j=n}^{+\infty}\left(V_j-\frac{1}{f(j)}\right)$$
Then there exists a constant $C=C_f>0$ such that for all large
enough $n\in \N$ and all $t\in\R^+$
\begin{eqnarray*}\Pr{A_{n}\;>\;\;t\sqrt{S_2(n)}}&\leq &
C e^{-t}\\
\Pr{A_{n}<-t\sqrt{S_2(n)}}&\leq & C e^{-t}\end{eqnarray*}
\end{lemma}
What \lemref{exp_lgdev} means to us is that $\sum_{j=a_1}^{+\infty}
X(1,j)$ can be though of as ``almost constant" in many calculations.
This will be put to use in all of our main proofs. In fact, the
following corollary will suffice.

\begin{corollary}\label{cor:exp_lgdev}Let $f$ and $\{V_j\}_{j\in\N}$ be as above, and define
$$B_n \equiv \sum_{j=n}^{+\infty} V_j\periodeq$$ Then $\Ex{B_n}=S_1(n)$ and there exists a $C'=C'_f>0$ a such
that
$$\forall n\in\N,\;\;\;\Pr{\left|\frac{B_n}{S_1(n)} - 1\right| >
\frac{C'}{n^{\frac{1}{4}}}}\leq C'
e^{-n^{\frac{1}{4}}}$$\end{corollary}

\begin{proof}[of \lemref{exp_lgdev}]  We will only prove the first
inequality, since the other proof is similar. The technique we
employ is fairly standard and is commonly used in other proofs of
Chernoff-type large deviation inequalities
\cite{AlonSpencer_Method}.

Fix any $0<s\leq f(n)/2=\min_{j\geq n}f(j)/2$ and notice that, by
the standard Bernstein's trick, the formulae in \secref{prelim},
the inequality ``$1+x\leq e^x$'', and some simple calculations
\begin{eqnarray*}\Pr{A_n\,>\,t\sqrt{S_2(n)}} &=& \Pr{e^{s\,A_n}>e^{s\,t\sqrt{S_2(n)}}}\\
& \leq &  e^{-s\,t\sqrt{S_2(n)}}\Ex{e^{\sum_{j\geq n}s \left(V_j - \frac{1}{f(j)}\right)}} \\
& = &  e^{-s\,t\sqrt{S_2(n)}}\prod_{j\geq n}\Ex{e^{s \left(V_j - \frac{1}{f(j)}\right)}}\\
& = &  e^{-s\,t\sqrt{S_2(n)}}\prod_{j\geq n}\frac{e^{-\frac{s}{f(j)}}}{1-\frac{s}{f(j)}} \\
& = &  e^{-s\,t\sqrt{S_2(n)}} \times \\ && \prod_{j\geq
n}e^{-\frac{s}{f(j)}}\left(1+
\frac{s}{f(j)} + \frac{s^2}{f(j)^2}\frac{1}{1-\frac{s}{f(j)}}\right) \\
&\leq & e^{-s\,t\sqrt{S_2(n)}}\prod_{j\geq n}\exp(2\frac{s^2}{f(j)^2})\\
\label{eq:est_lastexp}&=& \exp(2s^2\,S_2(n) -
s\,t\sqrt{S_2(n)})\end{eqnarray*} We now set
$$s\equiv \frac{1}{\sqrt{S_2(n)}}.$$
If we show that this choice is permissible (i.e. that
$1/\sqrt{S_2(n)}\leq f(n)/2$ for $n$ large), we will have finished
the proof. But notice that
$$\frac{1}{S_2(n)}\sim \frac{(2h(n)-1)f(n)^2}{n} \ll f(n)^2\commaeq$$
since $h(n)\ll n^{1/4}$ by assumption. We deduce that there is
some $n_0=n_0(f)$ such that for all $n\geq n_0$ $1/S_2(n)\leq
f(n)^2/4$, from which the lemma follows.\end{proof}
\begin{proof}[of \corref{exp_lgdev}] Recall that in this case
$0<B_n<+\infty$, since $B_n$ is positive and has finite expectation.
Hence, if $A_n$ is as in \lemref{exp_lgdev}, $B_n = A_n + S_1(n)$.
The Corollary follows from \lemref{exp_lgdev} by the choice of
$t\equiv n^{1/4}$ and recalling \lemref{est_moment}, which implies
that $$S_2(n)^{1/2}
=\bigoh{S_1(n)\sqrt{h(n)/n}}=\bigoh{n^{-1/4}S_1(n)} \mbox{ as
$h(n)\ll n^{1/4}$ by assumption.}$$\end{proof}

\section{Imbalanced start}\label{sec:tails_imbalanced}

We now start the discussion of the first of our main results. Recall
the notation $[n,\alpha]$ for the state of a two-bin balls-in-bins
process with feedback (cf. \secref{prelim}). Our interest in this
section will be in processes with two bins, started from
$[n,\alpha]$ with $n$ large, $\alpha\in[0,1/2)$ fixed, and $f$
valid. To state this theorem, we need a definition.
\begin{definition}\label{def:tails_hasmorethan}Let $\beta\in(0,1/2)$ and $N\in\N$ be given, and consider a balls-in-bins process with two bins. $\hasmorethan(\beta,N)$ is the event
that there are more than $\beta N$ balls in bin $1$ at the moment
when there is a total of $N$ balls in both bins.\end{definition}
\begin{theorem}\label{thm:tails_lgdev}Suppose that $f$ is a valid feedback function, so that in particular $\sum_{j\in\N}f(j)^{-1}<+\infty$ and monopoly has probability $1$ (cf. \propref{est_mon}).  Then for all $0<\alpha<\beta<1/2$, there exists a constant $\gamma>0$ depending only on
$\alpha$, $\beta$ and $f$ such that for all large enough $n\in \N$,
$$\Prp{[n,\alpha]}{\exists N\geq n\commaeq\;\hasmorethan(\beta,n)} \leq e^{-n^{\gamma}}\periodeq$$
In particular, if the initial conditions are $[n,\alpha]$ as above,
bin $2$ achieves monopoly with all-but-exponentially-small
probability.\end{theorem}

\begin{proof}[of \thmref{tails_lgdev}] Recall the definition of the exponential embedding: for $i=1,2$, $X(i,0)$ parameterizes the time until the first
arrival at bin $i$, while $X(i,j)$ (for $j>0$) parameterizes the
time between the $(j-1)$th and $j$th arrivals at bin $i$. We begin
by showing the following fact, which will also be useful in
\secref{tails_fraction}.

\begin{claim}\label{claim:tails_hasmorethan}Consider a balls-in-bins process with feedback
function $f$ and initial conditions $(x,y)$ with $x+y=n$. Let
$\beta\in(0,1)$ and $N\geq n$ be given. Then
\begin{equation}\label{eq:tails_hasmorethan}\hasmorethan(\beta,n) =
\left\{\sum_{j=x}^{\lceil \beta
N\rceil-1}X(1,j)<\sum_{\ell=y}^{N-\lceil \beta N\rceil-1}
X(2,\ell)\right\}\periodeq\end{equation}\end{claim}
\begin{proof}[of \claimref{tails_hasmorethan}]Let $A$ be the event in the RHS of \eqnref{tails_hasmorethan}. We begin by assuming that $\hasmorethan(\beta,N)$ occurs and show that this implies the occurrence of $A$. Consider the time $\tau_N$ when the total number of balls in the continuous time process reaches $N$. At that time, the number of balls in bin $1$ (respectively $2$) is larger than $\lceil \beta N\rceil$ (resp. smaller than $N-\lceil \beta N\rceil$), by assumption, so
$$\sum_{j=x}^{\lceil \beta
N\rceil-1}X(1,j)\leq \tau_N<\sum_{\ell=y}^{N-\lceil \beta N\rceil-1}
X(2,\ell)\commaeq$$ which implies the occurrence of $A$. Conversely,
assume that $A$ occurs, and let $\tau_N$ be as above. Then, because
$$\sum_{j=x}^{\lceil \beta N\rceil-1}X(1,j)<\sum_{\ell=y}^{N-\lceil \beta N\rceil-1}
X(2,\ell)\commaeq$$ the number of balls in bin $2$ at time
$\sum_{j=x}^{\lceil \beta N\rceil-1}X(1,j)$ is smaller than
$N-\lceil \beta N\rceil$, so that
$$\sum_{j=x}^{\lceil \beta N\rceil-1}X(1,j)<\tau_N\periodeq$$ This implies that at time $\tau_N$, the
number of balls at bin $1$ is at least $\lceil \beta N\rceil$, which
implies the occurrence of $\hasmorethan(\beta,N)$.\end{proof} We now
continue the proof of \thmref{tails_lgdev}. Given a number $D>0$
independent of $n$, consider the event $E_n$ where the four
conditions given below hold simultaneously:
\begin{eqnarray}\label{eq:tails_firstE}\left|\frac{\sum_{j=\lceil \alpha n\rceil}^{+\infty}X(1,j)}{S_1(\lceil \alpha n\rceil)}
-1\right| &\leq& D\,n^{-\frac{1}{4}}\commaeq\\
\forall N\geq n,\; \left|\frac{\sum_{j=\lceil \beta
N\rceil}^{+\infty}X(1,j)}{S_1(\lceil\beta N\rceil)} -1 \right|
&\leq& D\,N^{-\frac{1}{4}}\commaeq \\
\left|\frac{\sum_{\ell=n-\lceil \alpha
n\rceil}^{+\infty}X(2,\ell)}{S_1(n-\lceil \alpha n\rceil)} -1
\right| &\leq& D\, n^{-\frac{1}{4}} \commaeq\\
\label{eq:tails_lastE} \forall N\geq
n,\;\left|\frac{\sum_{\ell=N-\lceil \beta N\rceil
}^{+\infty}X(2,\ell)}{S_1(N-\lceil\beta N\rceil)} -1 \right| &\leq&
D\, N^{-\frac{1}{4}}\commaeq\end{eqnarray} The above conditions
correspond to the class of events covered by \corref{exp_lgdev}. For
instance, to get the first condition we may look at the
concentration of $B_{\lceil \alpha n\rceil}=\sum_{j=\lceil \alpha
n\rceil}^{+\infty}X(1,j)$. It follows that there exists a $D>0$
depending only on $C'=C'_f$ as in the Corollary, $\alpha$ and
$\beta$ for which:
$$\Pr{E_n}\geq 1 - 2C'e^{-n^{1/4}} - \sum_{N\geq n}2C'e^{-N^{1/4}}\geq 1 - e^{-n^\gamma}$$
with $\gamma>0$ depending only on $f,\alpha,\beta$.

From now on, our goal will be to show that:
\begin{equation}\label{eq:endfirstproof}\mbox{For all large enough $n\in\N$ and all $N\geq n$,}\; E_n\cap \hasmorethan(\beta,N)=\emptyset.\end{equation}
Notice that this implies the Theorem, as $$\Pr{\exists N\geq
n,\hasmorethan(\beta,N)}\leq \Pr{E_n^c}\leq e^{-n^\gamma}\mbox{ for
all large $n$.}$$

To establish \eqnref{endfirstproof} we note that there is nothing to
prove if $N-n<\lceil \beta N\rceil -\lceil\alpha n \rceil$: in that
case, there cannot be $\geq \lceil \beta N\rceil$ balls in bin $1$
after $N-n$ balls are added to the system. Hence we can assume that
\begin{equation}\label{eq:rangeofT}N-\lceil \beta N\rceil\geq n -\lceil\alpha n\rceil.\end{equation}
Inside $E_n$, we can use \eqnref{tails_firstE} --
\eqnref{tails_lastE}, \lemref{est_moment} and \lemref{est_bounded}
(to get rid of ceilings) to deduce:\begin{eqnarray*}\sum_{j=\lceil
\alpha n\rceil}^{\lceil \beta N\rceil -1}X(1,j) &=& \sum_{j=\lceil
\alpha n\rceil}^{+\infty}X(1,j) - \sum_{j=\lceil \beta N\rceil}^{+\infty}X(1,j)\\
&\geq& (1+\liloh{1})M_1(\alpha n) -(1+\liloh{1})M_1(\beta N)\\ &=&
(1+\liloh{1})\int_{\alpha n}^{+\infty}\frac{dx}{f(x)} -
(1+\liloh{1})\int_{\beta N}^{+\infty}\frac{dx}{f(x)},\end{eqnarray*}
where the $\liloh{1}$ terms converge to $0$ as $n\to +\infty$,
uniformly over $N$ satisfying \eqnref{rangeofT}. Similarly, we have
$$\sum_{j=n-\lceil \alpha n\rceil}^{N-\lceil \beta N\rceil}X(2,j)\leq
(1+\liloh{1})\int_{(1-\alpha) n}^{+\infty}\frac{dx}{f(x)} -
(1+\liloh{1})\int_{(1-\beta N)}^{+\infty}\frac{dx}{f(x)}.$$
Moreover, under initial conditions $[n,\alpha]$ one has
\begin{equation}\label{eq:hastooccur}\hasmorethan(\beta,N) = \left\{\sum_{j=\lceil \alpha n\rceil}^{\lceil \beta
N\rceil -1}X(1,j)<\sum_{\ell= n-\lceil \alpha n\rceil}^{N -\lceil
\beta N\rceil}X(2,\ell)\right\}.\end{equation} It follows that
$E_n\cap \hasmorethan(\beta,N)\neq \emptyset$ implies
$$(1+\liloh{1})\int_{\alpha n}^{+\infty}\frac{dx}{f(x)} -
(1+\liloh{1})\int_{\beta N}^{+\infty}\frac{dx}{f(x)}\leq
(1+\liloh{1})\int_{(1-\alpha) n}^{+\infty}\frac{dx}{f(x)} -
(1+\liloh{1})\int_{(1-\beta)N}^{+\infty}\frac{dx}{f(x)}.$$ This is
equivalent to: \begin{equation}\label{eq:cannotoccur}E_n\cap
\hasmorethan(\beta,N)\neq \emptyset\Rightarrow \int_{\alpha
n}^{(1-\alpha)n}\frac{dx}{f(x)} \leq (1+\liloh{1})\,\int_{\beta
N}^{(1-\beta)N}\frac{dx}{f(x)},\end{equation} since by
\lemref{est_moment} there is some $a<1$ such that
$$M_1(\alpha n)=\int_{\alpha
n}^{+\infty}\frac{dx}{f(x)}<a M_1((1-\alpha)
n)=a\,\int_{(1-\alpha)n}^{+\infty}\frac{dx}{f(x)}$$ and
$$M_1(\beta N)=\int_{\beta
N}^{+\infty}\frac{dx}{f(x)}<aM_1((1-\beta)n)=a\,\int_{(1-\beta)N}^{+\infty}\frac{dx}{f(x)}.$$

We will finish the proof by showing that the RHS of
\eqnref{cannotoccur} cannot hold for large $n$ and $N$ satisfying
\eqnref{rangeofT}. To this end, we employ estimate
\eqnref{est_fracgrowth} from the proof of \lemref{est_moment}, which
states that
\begin{equation}\label{eq:tails_fracgrowth}\forall y\geq x \;\; \frac{f(y)}{f(x)}\geq
\left(\frac{y}{x}\right)^c\commaeq\; \mbox{ where }c\equiv
\inf_{x'\geq x}h(x')\end{equation} We will only employ this
inequality for large $x<y$, which means we can assume $c>1$, because
$f$ is a valid feedback function. For any $N$ satisfying
\eqnref{rangeofT},
\begin{eqnarray*}\int_{\beta N}^{(1-\beta)N}\frac{dx}{f(x)} &\leq
&\left(\frac{n}{N}\right)^{c}\int_{\beta
N}^{(1-\beta)N}\frac{dx}{f\left(\frac{n}{N}\,x\right)} \\
&=& \left(\frac{n}{N}\right)^{c-1}\int_{\beta
n}^{(1-\beta)n}\frac{dy}{f\left(y\right)} \\
&=&
\left((1+\liloh{1})\frac{1-\beta}{1-\alpha}\right)^{c-1}\int_{\beta
n}^{(1-\beta)n}\frac{dy}{f\left(y\right)}\periodeq\end{eqnarray*}
Here, we employed \eqnref{tails_fracgrowth} for the second line, the
substitution $x\mapsto Ny/n$ in the third line, and
\eqnref{rangeofT} on the third, with (once again) a $\liloh{1}$ term
that is uniform over $n$. Since $c>1$ and $\beta>\alpha$ (hence
$1-\beta<1-\alpha$), there exists a constant $0<d<1$ such that for
all large enough $n$ and all $N$ satisfying \eqnref{rangeofT},
$$\int_{\beta N}^{(1-\beta)N}\frac{dx}{f(x)} \leq d\, \int_{\beta
n}^{(1-\beta)n}\frac{dy}{f\left(y\right)}\periodeq$$ To conclude,
note that $[\alpha,1-\alpha]\supset[\beta,1-\beta]$, hence the above
implies $$\int_{\beta N}^{(1-\beta)N}\frac{dx}{f(x)} \leq d\,
\int_{\alpha n}^{(1-\alpha)n}\frac{dx}{f\left(x\right)}\mbox{ with
$d<1$ constant},$$ which is incompatible with the RHS of
\eqnref{cannotoccur}. This finishes the proof.\end{proof}

\section{The number of balls in the losing bin}\label{sec:tails_losing}

We now prove the first of our two heavy-tail results. We first
recall the definition of $L$.

\begin{definition}Let $f$ be a feedback function satisfying the monopoly condition
$\sum_{j}f(j)^{-1}<+\infty$, so that in the corresponding
balls-in-bins process there almost surely is one bin that receives
all but finitely many balls. For a two-bin process with feedback
function $f$, the {\em losing number} $L$ is the (almost surely
finite) number of balls that go into the remaining
bin.\end{definition}

Our heavy tails result for $L$ is stated and proved below.

\begin{theorem}\label{thm:tails_losing}Let $f$ be a valid feedback function, in which case the monopoly condition $\sum_{j}f(j)^{-1}<+\infty$ is satisfied (cf. \propref{est_mon}). For any fixed initial conditions $(x,y)$, there exists a number $c>0$ (depending
only on $x$, $y$ and $f$) such that, as $n\to +\infty$,
$$\Prp{(x,y)}{L>n} \sim c\, S_1(n) \sim c\, \frac{n}{(h(n)-1)f(n)}\commaeq$$
where $h(n)$ is the characteristic exponent of $f$ (cf.
\defref{est_smooth} and \lemref{est_moment}).\end{theorem}
In the case $f(n)=n^p$, $p>1$, $S_1(n)\sim n^{1-p}/(p-1)$ for $n$
large, and as a consequence we have the following corollary.
\begin{corollary}\label{cor:tails_losing}For $f(n)\sim n^p$ with $p>1$ and $(x,y)\in\N^2$ fixed, as $n\to +\infty$
$$\Prp{(x,y)}{L>n} \sim \frac{c}{(p-1) \, n^{p-1}}\commaeq$$
with $c$ as above.\end{corollary}

\begin{proof}[of \thmref{tails_losing}]We will assume without loss of generality that $x\geq y$. First, notice
that for any $n>x$
\begin{eqnarray*}\label{eq:tail_Levent}\{L>n\}
 &=&\left\{\sum_{\ell=y}^{n}X(2,\ell) <
\sum_{j=x}^{\infty}X(1,j)< \sum_{\ell=y}^{\infty}X(2,\ell)
\right\} \\
& &\bigcup \left\{\sum_{j=x}^{n}X(1,j) <
\sum_{\ell=y}^{\infty}X(2,\ell)<\sum_{j=x}^{\infty}X(1,j)
\right\}\periodeq\end{eqnarray*} Indeed, $L>n$ if and only if bin
$1$ explodes first, but bin $2$ has at least $n$ balls when that
happens, or vice-versa. Define
\begin{eqnarray*}\label{eq:tails_Delta}\Delta_n &\equiv&
\sum_{j=x}^{n-1}X(1,j) - \sum_{\ell=y}^{n-1}X(2,\ell)\\
\label{eq:tails_A1} B^{(1)}_n &\equiv & \sum_{j=n}^{\infty}X(1,j)\\
\label{eq:tails_A2}B^{(2)}_n &\equiv &
\sum_{\ell=n}^{\infty}X(2,\ell)\end{eqnarray*} These random
variables are independent, and one can rewrite
$$\label{eq:tails_LDelta}\{L>n\} = \{0<\Delta_n<B^{(2)}_n\}\cup
\{0>\Delta_n>-B^{(1)}_n\}.$$

$B^{(1)}_n$ and $B^{(2)}_n$ have the same distribution and are
independent of $\Delta_n$, and the distribution of $\Delta_n$ has no
point masses. It follows that

\begin{eqnarray*}\nonumber \Prp{(x,y)}{L>n} &=&
\Pr{0<\Delta_n<B^{(2)}_n} +\Prp{(x,x)}{0>\Delta_n>-B^{(1)}_n}
\\ \nonumber &=&\Pr{0<|\Delta_n|<B_n} \\ \label{eq:tails_Lsym}  &=&\Pr{|\Delta_n|\leq B_n}\end{eqnarray*}
where we let $B_n\equiv B^{(1)}_n$ for simplicity.

We now apply the concentration result, \corref{exp_lgdev}, to
$B_n=\sum_{j=n}^{+\infty}V_j$ with $V_j=X(1,j)$. It follows that
there exists some $C'>0$ depending only on $f$ such that:
$$\label{eq:powerlaw_A_n} \Pr{\left|\frac{B_n}{S_1(n)} - 1\right|\geq
\frac{C'}{n^{1/4}}}\leq C'e^{-n^{1/4}}\periodeq$$ Plugging this into
the previous equation yields
\begin{eqnarray}\nonumber & \Pr{|\Delta_n|\leq (1-C'n^{-1/4})S_1(n)} - C'e^{-n^{1/4}} \\
\nonumber & \leq \Prp{(x,y)}{L>n}\leq \\ &\Pr{|\Delta_n|\leq
\label{tails_Lsym}(1+C'n^{-1/4})S_1(n)} +
C'e^{-n^{1/4}}\periodeq\end{eqnarray}

We will eventually show that there exist a sequence
$\{c_n\}_{n\in\N}$ and a constant $C>0$, both of which depend only
on $f$, $x$ and $y$, such that
\begin{equation}\label{eq:tails_delta2}|\Pr{|\Delta_n|\leq \eps} -
c_n \eps| \leq C \eps^3\commaeq\end{equation} and $c_n\to c$ for
some real-valued $c>0$. This inequality implies that, as $n\to
+\infty$, $$\label{eq:almostfinishes}|\Pr{|\Delta_n|\leq (1\pm
C'n^{-1/4})S_1(n)} - c_n\, S_1(n)|=\bigoh{n^{-1/4}S_1(n) +
S_1(n)^3}$$ Since we know that $e^{-n^{1/4}}\ll S_1(n)\ll 1$, this
implies the Theorem via \eqnref{tails_Lsym}.

To prove \eqnref{tails_delta2}, first let $\psi_n(\cdot)$ be the
characteristic function of $\Delta_n$.
\begin{eqnarray*}\psi_n(t) &=& \Ex{\exp(\I t\Delta_n)} \\ &=&
\Ex{\exp\left[\sum_{j=x}^{n-1}X(1,j)-\sum_{\ell=y}^{n-1} X(2,\ell) \right]} \\
& =& \left(\prod_{j=x}^{n-1} \frac{1}{1 - \frac{\I t}{f(j)}}\right)
\times \left(\prod_{j=y}^{n-1} \frac{1}{1 + \frac{\I
t}{f(\ell)}}\right) \\ &=& \left(\prod_{\ell=y}^{x-1} \frac{1}{1 +
\frac{\I t}{f(j)}}\right) \times
\left(\prod_{j=x}^{n-1}\frac{1}{1+\frac{t^2}{f(j)^2}}
\right)\end{eqnarray*} Clearly, for all $n\in\N$ $\Delta_n$ has a
distribution with no point masses, hence the inversion formula for
characteristic functions (a.k.a. Fourier inversion formula)
\cite{VaradhanBook} implies
\begin{eqnarray*}\forall -\infty<a\leq b<+\infty \;\; \Pr{a\leq \Delta_n\leq b} \\ = \frac{1}{2\pi} \lim_{T\nearrow +\infty} \int_{-T}^T \psi_n(t)\left(\frac{e^{-\I t a}-e^{-\I tb}}{\I t}\right)\,dt\periodeq\end{eqnarray*} In the present setting, we use this formula with
$b=-a=\eps$, to prove \eqnref{tails_delta2}, noting that, since
$\psi_n$ is integrable, we can dispense with the limit in $T$.
\begin{eqnarray*} \eps^{-1}\Pr{|\Delta_n|\leq \eps} &=&
\frac{1}{\pi}\lim_{T\nearrow +\infty} \int_{-T}^T
\psi_n(t)\left(\frac{e^{+\I\eps t}-e^{-\I \eps t}}{2\I \eps
t}\right)\, dt \\ &=&\frac{1}{\pi} \int_{-\infty}^{+\infty}
\psi_n(t)\frac{\sin(\eps t)}{\eps t}\, dt \end{eqnarray*} Now notice
that, quite crudely,
$$\forall s\in\R\commaeq\;|\frac{\sin(s)}{s}-1|\leq Cs^2$$
for some constant $C>0$. Applying this inequality with $s=\eps t$,
one obtains
\begin{eqnarray*}\left|\pi\eps^{-1}\Pr{|\Delta_n|\leq \eps}  - \int_{-\infty}^{+\infty} \psi_n(t)\,
dt\right|&\leq&\int_{-\infty}^{+\infty}
|\psi_n(t)|\left|\frac{\sin(\eps t)}{\eps t} - 1\right|\, dt \\
&\leq &C\eps^2\, \int_{-\infty}^{+\infty} |t^2\,\psi_n(t)|\,
dt\periodeq\end{eqnarray*} For $n=x+3$, $|t^2\psi_n(t)|$ is of order
$1/t^2$, so the above integral converges; for $n>x+3$, the integrand
is even smaller. Hence, we can guarantee that, for a possibly larger
$C>0$, $$\left|\Pr{|\Delta_n|\leq \eps}  -
\left(\frac{1}{\pi}\int_{-\infty}^{+\infty} \psi_n(t)\, dt
\right)\eps \right| \leq C \eps^3\mbox{ uniformly over
$n$}\periodeq$$ Moreover, since $|\psi_{n+1}(t)|\leq |\psi_{n}(t)|$
for all $t$, the Dominated Convergence Theorem implies that
\begin{equation}\label{eq:tails_defcn} c_n\equiv \frac{1}{\pi}\int_{-\infty}^{+\infty} \psi_n(t)\,
dt\to c \equiv
\frac{1}{\pi}\int_{-\infty}^{+\infty}\left(\prod_{\ell=y}^{x-1}
\frac{1}{1 + \frac{\I t}{f(j)}}\right) \times
\left(\prod_{j=x}^{n-1}\frac{1}{1+\frac{t^2}{f(j)^2}} \right)\,
dt\periodeq\end{equation}

Our last step is to prove that $c=\lim_n c_n>0$. To see this,
consider first the case $x=y$. In this case, the formula for $c$ in
\eqnref{tails_defcn} becomes
$$\frac{1}{\pi}\int_{-\infty}^{+\infty}\left(\prod_{j=x}^{+\infty}\frac{1}{1+\frac{t^2}{f(j)^2}} \right)\,
dt\periodeq$$ The product in the integrand converges to a positive
limit for all $t\in\R$, since $\sum f(j)^{-1}<+\infty$; hence, the
value of the integral is positive, and we are done in this case.\\

We now consider the case $y<x$. Clearly,
$c=\lim_{n\to+\infty}S_1(n)^{-1}\Pr{L>n}$ is a real number.
Moreover, notice that there is a positive probability $\alpha$
that in the process started from $(x,y)$, bin $2$ receives the
first $x-y$ balls, thus evolving to state $(x,x)$. Conditioned on that happening, the probability
of $L>n$ is $\Prp{(x,x)}{L>n} \sim c'S_1(n)$ for some $c'>0$ (as
shown above). But then
\begin{eqnarray*}c&=&\lim_{n\to+\infty}\frac{\Prp{(x,y)}{L>n}}{S_1(n)} \\
&\geq& \Prp{(x,y)}{\mbox{bin $2$ gets first $x-y$ balls}}\\ \nonumber & & \times  \lim_{n\to+\infty}\frac{\Prp{(x,y)}{L>n\mid \mbox{bin $2$ gets first $x-y$ balls}}}{S_1(n)}\\
&=&\alpha \lim_{n\to+\infty}\frac{\Prp{(x,x)}{L>n}}{S_1(n)} =\alpha
c'>0\end{eqnarray*} which proves that $c$ is positive even when $x$
and $y$ differ, thus finishing the proof.\end{proof}

\section{The time until imbalance}\label{sec:tails_fraction}

The strategy used to prove \thmref{tails_losing} was very simple.
After the event $\{L>n\}$ was written down in terms of the
exponential embedding random variables, $B_n$ was approximated by
its expectation, and the distribution of $\Delta_n$ near the origin
via Fourier transform techniques. As we shall see below, our last
theorem in this section has a similar proof. First, we need a
definition.
\begin{definition}Given a number $n\in\N$, a number $0<\alpha<1/2$ and a balls-in-bins process with two bins, the event $\loserhasmorethan(\alpha,n)$ holds if at the time the number of balls in the system reaches $n$, the number of balls in the losing bin is at least $\alpha n$.\end{definition}
Our second heavy-tails result can now be properly stated.

\begin{theorem}\label{thm:tails_fraction}Let $f$ be a valid feedback function, in which case the monopoly condition $\sum_{j}f(j)^{-1}<+\infty$ is satisfied (cf. \propref{est_mon}). Then for all fixed
initial conditions $(x,y)\in \N^2$ and $0<\alpha<1/2$ there is a
constant $c>0$ depending only on $f$, $x$ and $y$ such that
$$\Prp{(x,y)}{\loserhasmorethan(\alpha,n)} \sim c [S_1(\lceil\alpha n \rceil,n-\lceil \alpha n \rceil)]\mbox{ as }n\to+\infty\periodeq$$
Moreover, we can take $c$ to be the same constant (depending on
$f$, $x$ and $y$) that appears in the proof of
\thmref{tails_losing}. \ignore{ Finally, if $f$ is * and satisfies
the * condition,
$$\Prp{(x,y)}{\exists N\geq n\,\, \loserhasmorethan(\alpha,n)} \sim c [S_1(\lceil\alpha n \rceil )-S_1(n-\lceil \alpha n \rceil)] \mbox{ as
}n\to+\infty\periodeq$$}\end{theorem}

In the case $f(x)=x^p$ for $p>1$, the estimate of $S_1(n)$ by an
integral implies the following corollary.

\begin{corollary}\label{cor:tails_fraction}For $f(n)=n^p$ with $p>1$ and $(x,y)\in\N^2$, $0<\alpha<1/2$ fixed, as $n\to +\infty$
$$\Prp{(x,y)}{\loserhasmorethan(\alpha,n)} \sim \frac{c\left(\alpha^{1-p}-(1-\alpha)^{1-p}\right)}{(p-1) \, n^{p-1}}\commaeq$$
with $c$ as above.\end{corollary}

\begin{proof}[of \thmref{tails_fraction}] We follow the same outline as in the proof of \thmref{tails_losing}. We
will again assume that $x\geq y$, and write the event under
consideration in terms of the random variables in the definition of
the exponential embedding. To do that, first notice that
$\loserhasmorethan(\alpha,n)$ occurs if and only if both bins have
at least $\alpha n$ balls at the time when the total number of balls
in the system is $n$. Indeed, if one of the bins has less than
$\alpha n<n/2$ balls, then this must necessarily be the losing bin,
and $\loserhasmorethan(\alpha,n)$ cannot occur in this case.
Conversely, if both bins have at least $\alpha n$ balls, then in
particular the losing bin has $\geq \alpha n$ balls, and
$\loserhasmorethan(\alpha,n)$ occurs.\\

The event that bin $1$ has at least $\alpha n$ balls when $n$ balls
are present in the system is precisely the event
$\hasmorethan(\alpha,n)$ defined in
\defref{tails_hasmorethan}, which was shown in
\claimref{tails_hasmorethan} to be equal to
$$\hasmorethan(\alpha,n)=\left\{\sum_{j=x}^{\lceil \alpha n\rceil
-1}X(1,j) < \sum_{\ell=y}^{n-\lceil \alpha n\rceil
-1}X(2,\ell)\right\}\periodeq$$ Similarly, the event than bin $2$
has at least $\alpha n$ balls when there are $n$ balls in the system
is precisely $\hasmorethan(\alpha,n)$ with the roles of the two bins
reversed, and can thus be written down as
$$\left\{\sum_{\ell=y}^{\lceil \alpha n\rceil -1}X(2,\ell) <
\sum_{j=x}^{n-\lceil \alpha n\rceil -1}X(1,j)\right\}\periodeq$$ We
conclude that
\begin{eqnarray}\nonumber & & \loserhasmorethan(\alpha,n)\\
\nonumber &=& \left\{\sum_{j=x}^{\lceil \alpha n\rceil -1}X(1,j) <
\sum_{\ell=y}^{n-\lceil \alpha n\rceil -1}X(2,\ell)\right\} \\
\nonumber & & \bigcap \left\{\sum_{\ell=y}^{\lceil \alpha n\rceil
-1}X(2,\ell) < \sum_{j=x}^{n-\lceil \alpha n\rceil -1}X(1,j)\right\}
\\ \nonumber &= & \left\{\sum_{j=x}^{\lceil \alpha n\rceil -1}X(1,j) -
\sum_{j=y}^{\lceil \alpha n\rceil -1}X(2,\ell) < \sum_{\ell=\lceil
\alpha n\rceil}^{n-\lceil \alpha n\rceil -1}X(2,\ell)\right\} \\
\label{eq:tails_losefhasmorethan}& & \bigcap
\left\{\sum_{\ell=y}^{\lceil \alpha n\rceil -1}X(2,\ell) -
\sum_{j=x}^{\lceil \alpha n\rceil -1}X(1,j) < \sum_{j=\lceil \alpha
n\rceil}^{n-\lceil \alpha n\rceil
-1}X(1,j)\right\}\commaeq\end{eqnarray} where for the last equality
we used the fact that $\alpha<1/2$, which guarantees $\lceil \alpha
n\rceil <n/2$ for all large $n$, to ensure that the sums from
$\lceil\alpha n\rceil$ to $n-\lceil\alpha n\rceil$ are non-empty. If
we define:
\begin{eqnarray*}\Sigma_{n}&= &\sum_{j=x}^{\lceil \alpha n\rceil -1}X(1,j) -
\sum_{j=y}^{\lceil \alpha n\rceil -1}X(2,\ell)\commaeq \\
E^{(1)}_{n} & = &\sum_{j=\lceil \alpha n\rceil}^{n-\lceil
\alpha n\rceil -1}X(1,j)\commaeq \\
E^{(2)}_{n} & = &\sum_{\ell=\lceil \alpha n\rceil}^{n-\lceil \alpha
n\rceil -1}X(2,\ell)\commaeq\end{eqnarray*} we can rewrite
$\loserhasmorethan(\alpha,n)$ as
\begin{equation}\loserhasmorethan(\alpha,n) = \{-E^{(1)}_n < \Sigma_n <
E^{(2)}_n\}\periodeq\end{equation}The random variable $\Sigma_n$
equals $\Delta_{\lceil\alpha n\rceil}$, as defined in the proof of
\thmref{tails_losing}, and the $E^{(i)}_n$'s are akin to the
$B^{(i)}_n$'s in that proof. Similarly to that proof, we note that
$E^{(1)}_n$ and $E^{(2)}_n$ are independent, identically distributed
and independent from $\Sigma_n$. Since $\Sigma_n$ has no
point-masses in its distribution,
\begin{eqnarray*}& & \Prp{(x,y)}{\loserhasmorethan(\alpha,n)} \\ &=& 1 - \Pr{-E^{(1)}_n > \Sigma_n} -\Pr{\Sigma_n
>E^{(2)}_n} \\ &=& 1 - \Pr{-E_n > \Sigma_n} -\Pr{\Sigma_n
>E_n} \\ &=& \Pr{|\Sigma_n|<E_n}\commaeq\end{eqnarray*}
where $E_n=E^{(1)}_n$. Now notice that:
\begin{enumerate}
\item $E_n$ concentrates around its mean. Indeed,
$$E_n=\left[\sum_{j=\lceil \alpha n\rceil}^{+\infty}X(1,j)\right] - \left[\sum_{j=n-\lceil \alpha
n\rceil}^{+\infty}X(1,j)\right]\periodeq$$ Applying
\corref{exp_lgdev} to each bracketed term and noticing that (by
\lemref{est_moment})
$$S_1(\lceil \alpha n\rceil,n-\lceil \alpha n\rceil)=\ohmega{S_1(\lceil \alpha n\rceil)}\commaeq$$
we conclude that there exists a $D'>0$ depending only on $f$ and
$\alpha$ such that $$\Pr{\left|E_n - S_1(\lceil \alpha
n\rceil,n-\lceil \alpha n\rceil)\right| \geq
\frac{D'}{n^{1/4}}S_1(\lceil \alpha n\rceil)}\leq D'e^{-n^{1/4}}.$$
\item The estimates on $\Delta_n$ in \eqnref{tails_delta2}
imply that there exists a constant depending only on $f$ such that
for all $\eps>0$ and all $n$ large enough
$$\left|\Pr{|\Sigma_n|<\eps} - c_{\lceil\alpha n\rceil}\eps\right|
\leq C\eps^3$$for the same sequence $\{c_m\}_{m\in\N}$ converging to
$c>0$ appearing in the proof of \thmref{tails_losing}.
\end{enumerate} Putting those results together we conclude that
\begin{eqnarray*}&\left|\Pr{|\Sigma_n|\leq E_n} - c_{\lceil\alpha
n\rceil}[S_1(\lceil \alpha n\rceil,n-\lceil \alpha n\rceil)]\right|
\\ = & \bigoh{\frac{1}{n^{1/4}}[S_1(\lceil \alpha n\rceil) +
S_1(n-\lceil \alpha n\rceil)]+ [S_1(\lceil \alpha n\rceil) -
S_1(n-\lceil \alpha n\rceil)]^3 + e^{-n^{1/4}}}
\\ = & \liloh{[S_1(\lceil \alpha n\rceil)}\end{eqnarray*} which implies that
$$\Prp{(x,y)}{\loserhasmorethan(\alpha,n)}= \Pr{|\Delta_n|\leq
C_n}\sim c\,[S_1(\lceil \alpha n\rceil,n-\lceil \alpha
n\rceil)]\periodeq$$This is precisely the desired result.\end{proof}

\section{Extensions, related results and open problems}\label{sec:last}

\begin{itemize}
\item The proof of \thmref{tails_losing} generalizes directly to
the following statement.
\begin{theorem}\label{thm:losergeneral}Let $q:\N\to \N$ be a
function and $f$ be a valid feedback function such that
$$\label{eq:conditionforlast}S^2_1\left(\left\lceil\frac{n-q(n)}{2}\right\rceil,\left\lfloor\frac{n+q(n)}{2}\right\rfloor\right)
\gg n^{\gamma
}S_2\left(\left\lceil\frac{n-q(n)}{2}\right\rceil,\left\lfloor\frac{n+q(n)}{2}\right\rfloor\right)$$as
$n\to +\infty$, for some constant $\gamma>0$ depending only on $f$.
Then for any fixed $(x,y)\in\N^2$ there exists a constant $c>0$
depending only on $f$, $x$ and $y$ such that for $n\gg 1$
$$\Prp{(x,y)}{|I_1(n-(x+y))-I_2(n-(x+y))|\leq q(n)} \sim c\, S_1\left(\left\lceil\frac{n-q(n)}{2}\right\rceil,\left\lfloor\frac{n+q(n)}{2}\right\rfloor\right).$$
\end{theorem}
\thmref{tails_losing} is a special case of this result when
$q(n)=\alpha n$, and the r\^ole of \eqnref{conditionforlast} is to
show that $\sum_{j=(n-q(n))/2}^{(n+q(n))/2}X(i,1)$ concentrates
around its mean (cf. \lemref{exp_lgdev}). We omit the proof, but
note the following corollary (take $q(n)=\lambda \sqrt{n}$.)
\begin{corollary}Assume $f(n)=n^p$ for $p>1$. Then the probability that the losing bin has at least $(n-\lambda\sqrt{n})/2$ balls at the time when the total number of balls is $n$ is asymptotic to
$$(\mbox{const.})\times n^{1/2-p} \;\;\;(n\gg 1),$$
with a constant that depends on the initial conditions, $p$ and
$\lambda$.\end{corollary} We mention this corollary because it
relates to a result of \cite{MitzenmacherOS04}, where it is shown
that for a large initial number of balls $n$, the probability of
monopoly by bin $1$ has non-trivial behavior when
$|I_1(0)-I_2(0)|=\theita{\sqrt{n}}$.

\item We do not know any result more precise that
\thmref{tails_lgdev} about the behavior of
$\Prp{[t,\alpha]}{\hasmorethan(\beta,T)}$, but a related question
has been addressed. Let $f(x)=x^p$ and recall that $\mon_1$ is the
event that bin $1$ achieves monopoly. Assume we start with bin $1$
losing, from state $[t,\alpha]$, $0<\alpha<1/2$. We prove in
\cite{avoid} that
$$\Prp{[t,\alpha]}{\mon_1} = \exp(c_p(\alpha)t+o(t))$$
for some negative, smooth function $c_p(\cdot)<0$. Moreover, we
show that conditioned on $\mon_1$, the fraction of balls in the
first bin approximately follows the solution to a deterministic
ODE.
\end{itemize}

\bibliography{mybibtek}
\bibliographystyle{plain}

\end{document}